\documentclass{gtpart}
\usepackage[all]{xy}

\title{Homotopy theory of presheaves of $\Gamma$-spaces}

\author[H\,S Bergsaker]{H{\aa}kon Schad Bergsaker}
\givenname{H{\aa}kon}
\surname{Bergsaker}
\address{Department of Mathematics\\University of Oslo\\\newline Oslo\\Norway}
\email{hakonsb@math.uio.no}

\keyword{model categories}
\keyword{homotopy theory}
\keyword{simplicial presheaves}
\keyword{infinite loop spaces}

\subject{primary}{msc2000}{55P47}
\subject{secondary}{msc2000}{55P42}
\subject{secondary}{msc2000}{55P43}
\subject{secondary}{msc2000}{55P48}

\newtheorem{thm}{Theorem}[section]
\newtheorem{lem}[thm]{Lemma}
\newtheorem{cor}[thm]{Corollary}
\newtheorem{prop}[thm]{Proposition}
\newtheorem{hyp}[thm]{Hypothesis}

\theoremstyle{definition}
\newtheorem{defn}[thm]{Definition}   
\newtheorem*{rem}{Remark}            
\newtheorem*{notn}{Notation}
\newtheorem*{ack}{Acknowledgements}

\newcommand{\set}{\operatorname{Set}_*}
\newcommand{\fset}{\operatorname{fSet}_*}
\newcommand{\sset}{\mathcal{S}_*}
\newcommand{\gsset}{\Gamma\sset}
\newcommand{\site}{\mathcal{C}}

\newcommand{\fun}{\operatorname{Fun}}
\newcommand{\spc}{\operatorname{Spc}}
\newcommand{\sspc}{\operatorname{s}\spc}
\newcommand{\gspc}{\Gamma\spc}
\newcommand{\gmod}{\Gamma\operatorname{Mod}}
\newcommand{\galg}{\Gamma\operatorname{Alg}}
\newcommand{\sabpre}{\operatorname{sAbPre}}
\newcommand{\smodpre}{\operatorname{sModPre_A}}
\newcommand{\salgpre}{\operatorname{sAlgPre_B}}
\newcommand{\cat}{\mathcal{M}}

\newcommand{\vcat}{\mathcal{V}}
\newcommand{\op}{\operatorname{op}}

\newcommand{\gae}{\Phi(E)}
\newcommand{\colim}{\operatornamewithlimits{colim}}
\newcommand{\shom}{\operatorname{Map}}
\newcommand{\vhom}{\operatorname{\mathcal{V}Hom}}
\newcommand{\phom}{\operatorname{SpcHom}}
\newcommand{\ihom}{\underline{\operatorname{Hom}}}

\newcommand{\ev}{\operatorname{Ev}_n}
\newcommand{\lev}{\operatorname{L}_n}
\newcommand{\levm}{\operatorname{L}_m}
\newcommand{\levmn}{\operatorname{L}_{mn}}
\newcommand{\levnn}{\operatorname{L}_{n_1n_2}}
\newcommand{\spt}{\operatorname{Spt}}
\newcommand{\sptsym}{\operatorname{Spt}^\Sigma}
\newcommand{\ho}{\operatorname{Ho}}
\newcommand{\Sp}{\operatorname{Sp}}
\newcommand{\s}{\mathbb{S}}

\newcommand{\wF}{\bar F}

\numberwithin{equation}{section}

\begin{document}

\begin{abstract}
We consider the category of presheaves of $\Gamma$-spaces, or equivalently, of $\Gamma$-objects in simplicial presheaves.  Our main result is the construction of stable model structures on this category parametrised by local model structures on simplicial presheaves.  If a local model structure on simplicial presheaves is monoidal, the
corresponding stable model structure on presheaves of $\Gamma$-spaces is monoidal and satisfies the monoid axiom.
This allows us to lift the stable model structures to categories of algebras and modules over commutative algebras.
\end{abstract}

\begin{asciiabstract}
We consider the category of presheaves of Gamma-spaces, or equivalently, of 
Gamma-objects in simplicial presheaves.  Our main result is the construction of 
stable model structures on this category parametrised by local model structures
on simplicial presheaves.  If a local model structure on simplicial presheaves 
is monoidal, the corresponding stable model structure on presheaves of 
Gamma-spaces is monoidal and satisfies the monoid axiom.  This allows us to 
lift the stable model structures to categories of algebras and modules over 
commutative algebras.
\end{asciiabstract}

\maketitle


\section*{Introduction}

In his paper \cite{segal} Segal introduced $\Gamma$-spaces as a way to describe commutative monoids up to homotopy, and showed that they give rise to infinite loop spaces.
Segal's original definition of a $\Gamma$-space, as a functor from the category of finite sets to spaces satisfying
certain conditions, is what is now called a special $\Gamma$-space.  In \cite{bousfield-friedlander} Bousfield and Friedlander considered the category of all based functors from finite sets to simplicial
sets; and in particular constructed a stable model structure on it, in which the fibrant objects are given by
the very special $\Gamma$-spaces, and
the weak equivalences are the stable equivalences of the associated spectra.  As a consequence they show that the homotopy category of this model category is equivalent to the homotopy category of connective spectra.

Lydakis introduced a smash product for $\Gamma$-spaces in \cite{lydakis}, making the category of $\Gamma$-spaces
into a symmetric monoidal category.  This smash product is compatible with the smash product of spectra after passage to the respective homotopy categories, thus making the category of $\Gamma$-spaces a convenient category
for modeling connective spectra on a point set level.  In \cite{schwede}, Schwede introduced a different model
structure for $\Gamma$-spaces, Quillen equivalent to the one considered by Bousfield and Friedlander.  This model
structure satisfies the monoid axiom, an axiom first formulated by Schwede and Shipley in \cite{schwede-shipley}, 
which implies the existence
of model structures on the categories of monoids and modules of $\Gamma$-spaces.

The main result of this paper is the construction of stable model structures on the category of presheaves of $\Gamma$-spaces, 
or equivalently, of $\Gamma$-objects in simplicial presheaves over an arbitrary small Grothendieck site.  
There are several model structures on simplicial presheaves, and we are focusing on the ones with local weak 
equivalences (\ref{defn-localweak}) as weak equivalences.  We carry out the arguments without assuming any 
particular choice of model structure on simplicial presheaves, but have to impose a cofibrancy condition 
on the domains of the generating sets (Hypothesis \ref{spc-hyp}).
When the site consists of one morphism only, our model structure will specialize to the one in \cite{schwede}.

The following theorem states the main results appearing as Theorem
\ref{thm-stable-model}, Proposition \ref{prop-stable-monoidal} and Proposition \ref{prop-stable-monoid} in the 
main body of the paper.

\begin{thm}
\label{main-thm}
\begin{enumerate}
\item
Let $\site$ be a small Grothendieck site and let $\spc$ be the category of simplicial presheaves given a
model structure according to Hypothesis \ref{spc-hyp}.  Let $\gspc$ denote the category of based functors
$\Gamma \to \spc$, where $\Gamma$ is the category of finite ordinals.
There is a cofibrantly generated left proper model structure on the category $\gspc$ with stable equivalences 
(Definition \ref{defn-stable}) as weak equivalences.  The fibrant objects in 
this model structure coincides with the very special (Definition \ref{defn-special}) $\Gamma$-spaces.
\item
If the category $\spc$ is a monoidal model category, then the stable model structure on $\gspc$ is monoidal and
satisfies the monoid axiom.  Consequently, the category of algebra objects in $\gspc$, and the category of module objects over a commutative algebra in $\gspc$, inherits model structures from $\gspc$ by the results of
\cite{schwede-shipley}.
\end{enumerate}
\end{thm}

As a part of the construction, we
compare our $\Gamma$-spaces to presheaves of spectra, and also show that the homotopy category of (presheaves of) 
$\Gamma$-spaces is equivalent to the homotopy category of connective (presheaves of) spectra.
This equivalence is induced by a left Quillen functor from $\Gamma$-spaces to spectra which maps the very special
$\Gamma$-spaces to $\Omega$-spectra, thereby producing infinite loop objects in the category of simplicial presheaves.

As an application of the last part of Theorem \ref{main-thm} we construct an Eilenberg-Mac Lane functor $H$ from
presheaves of simplicial abelian groups to $\Gamma$-spaces and show that it is a Quillen equivalence between 
the categories of presheaves of simplicial abelian groups and the category of $H\Z$-modules.
Corresponding results for presheaves of simplicial rings, and presheaves of simplicial modules over presheaves of commutative simplicial rings are also included.  

Here is a quick outline of the paper.  In Section 1 we recall some basic theory of simplicial presheaves, in
particular the relevant model structures.  Section 2 introduces the category of $\Gamma$-spaces, and in
Section 3 we establish the pointwise model structure on this category.  We apply Bousfield localization to this 
model structure 
in Section 4 to obtain the stable model structure on $\Gamma$-spaces, and compare its homotopy category to
the homotopy category of connective presheaves of spectra.  In Section 5 the stable model structure is lifted to the
categories of modules and algebras over a (commutative) $\Gamma$-ring.  Here we also obtain a Quillen equivalence
between presheaves of simplicial modules over a presheaf of simplicial rings and modules over a $\Gamma$-ring.
To this end we first construct a model structure on presheaves of simplicial modules.  Similarly for algebras.

We assume familiarity with the theory of model categories, as described in e.g. Goerss and Jardine \cite{gj}, 
Hirschhorn \cite{hirschhorn} or Hovey \cite{hovey-model}.  
Some knowledge of classical $\Gamma$-spaces and simplicial presheaves is also assumed, but
we recall what we need about simplicial presheaves in the first section.
To prove the main theorem we make use of enriched left Bousfield localization as described in Barwick \cite{barwick}.
A quick review of this theory, together with some notes on bisimplicial presheaves, is located in an appendix.

\begin{ack}
This paper is part of my PhD thesis done at the University of Oslo, the topic was suggested to me by Paul Arne {\O}stv{\ae}r.  I want to thank Clark Barwick, John Rognes and Paul Arne {\O}stv{\ae}r for clarifying conversations
regarding this paper and model categories in general.
\end{ack}

\begin{notn}
In this paper, we use $\cat(X,Y)$ to denote the set of morphisms between $X$ and $Y$ in the category $\cat$, while
$\shom(X,Y)$ and $\ihom(X,Y)$ will denote respectively simplicial function complex and internal hom.  
More generally, when $\cat$ is enriched in a category $\vcat$, the enriched hom objects will be denoted $\vhom(X,Y)$.
When more than one category is under consideration, these objects will often be subscripted by the categories.
\end{notn}

\section{Preliminaries on spaces}

In this section we recall some facts about simplicial presheaves.
Let $\sset$ be the category of pointed simplicial sets. 
Fix a small site $\site$, i.e., a small category $\site$ with a Grothendieck topology.  The functor category
$\fun(\site^{op},\sset)$, which we denote $\spc$, is the category of pointed
simplicial presheaves on $\site$.  As the notation suggests, we will call the objects in this category
``spaces''.

Each $U \in \site$ represents a discrete simplicial presheaf $\site(-,U)$, and we will write $U$ for
this space.  Also, a simplicial set $K$ defines a constant simplicial presheaf and we 
will use $K$ to denote this space.

The category of spaces is closed symmetric monoidal, with monoidal product $\wedge$
defined sectionwise by
$$
(X \wedge Y)(U) = X(U) \wedge Y(U)
$$
for all $U \in \site$.
Here we are using $\wedge$ to denote both the monoidal product of spaces $X$ and
$Y$ and the smash product of based simplicial sets.  Let $K$ be a based simplicial set.  
Simplicial tensor $K \wedge -$ and cotensor $(-)^K$ are defined as
\begin{eqnarray*}
(K \wedge X)(U) = K \wedge X(U) \\
X^K(U) = X(U)^K
\end{eqnarray*}
for each $U \in \site$.

The simplicial function complex 
$\shom(X,Y)$ of two spaces $X$ and $Y$ is defined in simplicial degree $n$ to be
$$
\shom(X,Y)_n = \spc(X \wedge \Delta^n_+,Y) \,,
$$
with face and degeneracy maps induced from $\Delta^n_+$.
There is also an internal hom-object $\ihom(X,Y)$ of spaces defined sectionwise
by
$$
\ihom(X,Y)(U) = \shom(X|U,Y|U) \,,
$$
where $X|U$ means $X$ restricted to the local site $\site \downarrow U$.

We define homotopy groups of a space $X$ as follows.  First, let 
$$
L^2 \colon \text{Pre}(\site) \to \text{Shv}(\site) 
$$
be the associated sheaf functor from the category of presheaves to
the category of sheaves, which is left adjoint to the inclusion functor. 
Let $\pi^p_0(X)$ be the presheaf
$U \mapsto \pi_0(X(U))$; the sheaf of path components is the assoctiated sheaf $\pi_0(X) = L^2\pi^p_0(X)$.
For $n\geq 1$, each $U \in \site$ and $0$-simplex $x \in X(U)$, define the presheaf
$\pi_n^p(X,x)$ on $\site \downarrow U$ as 
$$ 
\pi_n^p(X,x)(V) = \pi_n(|X(V)|,x|V) \,, 
$$ 
where $|-|$ denotes geometric realization of simplicial sets and $x|V$ denotes the 
restriction of $x$ along $X(U) \to X(V)$.
The sheaf $\pi_n(X,x) = L^2\pi_n^p(X,x)$ is the sheaf of homotopy groups of $X$ 
over $U$ with basepoint $x$. 

\begin{defn}
\label{defn-localweak}
A morphism $f \colon X \to Y$ of spaces is a \emph{local weak equivalence} if the
induced map of sheaves $\pi_0(X) \to \pi_0(Y)$ is a bijection, and the induced maps
$$
\pi_n(X,x) \to \pi_n(Y,f(x))
$$
are isomorphisms for all $n\geq 1, U \in \site$, $x \in X(U)_0$.  It is a \emph{sectionwise
weak equivalence} if $f(U) \colon X(U) \to Y(U)$ is a weak equivalence of
simplicial sets for each $U \in \site$, and a sectionwise equivalence is in particular a local weak equivalence.
Sectionwise cofibrations and fibrations are defined similarly.  
\end{defn}

There are several known model structures on $\spc$.  We will only consider
model structures on $\spc$ in which the weak equivalences are given by the
local weak equivalences of spaces.

\begin{thm}[Jardine \cite{jardine-spre}]
\label{spc-inj-model}
There is a cofibrantly generated proper simplicial model structure on $\spc$
with sectionwise cofibrations (i.e., monomorphisms) as cofibrations and local weak equivalences as weak
equivalences.  This is the local injective model structure on $\spc$.
\end{thm}

To formulate the next theorem, let us define a \emph{projective cofibration} of spaces to be a map that
has the left lifting property with respect to maps that are both sectionwise fibrations and sectionwise
weak equivalences.

\begin{thm}[Blander \cite{blander-spre}]
\label{spc-proj-model}
There is a cofibrantly generated proper simplicial model structure on $\spc$
with cofibrations the projective cofibrations of spaces, and local weak equivalences as weak
equivalences.  This is the local projective model structure on $\spc$.
\end{thm}

Each projective cofibration $i \colon A \to B$ can be factored as a monomorphism $j \colon A \to C$ followed by a 
local injective trivial fibration $p \colon C \to B$.  Since $p$ is also a local projective trivial fibration there 
is a lift in the diagram
$$
\xymatrix{
A \ar[d]_i \ar[r]^j & C \ar[d]_p \\
B \ar[r] & B \,,
}
$$
from which we see that $i$ is a retraction of $j$, hence $i$ is a monomorphism.
This shows that the class of projective cofibrations is contained in the class of local
injective cofibrations.  In fact any set $I$ of monomorphisms
containing the set of generating projective cofibrations determines a local
model structure on $\spc$.

\begin{thm}[Jardine \cite{jardine-inter}]
\label{spc-inter-model}
Let $I$ be a set of monomorphisms conitaining the set of generating projective cofibrations.
There is a cofibrantly generated proper simplicial model structure on $\spc$ with
$I$ as the set of generating cofibrations and local weak equivalences as weak
equivalences.   
\end{thm}

An example of an intermediate model structure which differs from the local injective and
local projective ones is the flasque model structure constructed by Isaksen \cite{isaksen}.

\begin{prop}
\label{spc-smash}
If $f \colon X \to Y$ is a local weak equivalence of spaces and $Z$ is a space, then the induced
map $f \wedge 1 \colon X \wedge Z \to Y \wedge Z$ is a local weak equivalence.
\end{prop}
\begin{proof}
This is stated in \cite[2.46]{jardine-gen}.
\end{proof}

\begin{prop}
\label{spc-loc}
If $f \colon X \to Y$ is a local weak equivalence of spaces, where $X$ and $Y$ are fibrant in any of
the model structures constructed in Theorem \ref{spc-inter-model}, then $f$ is a sectionwise weak
equivalence.
\end{prop}
\begin{proof}
By \cite[1.3]{blander-spre} the local projective model structure on spaces is a Bousfield localization of 
the projective model structure consisting of the projective cofibrations and sectionwise fibrations and weak
equivalences, so in this case the result follows from general properties of Bousfield localizations.
But a space $X$ which is fibrant in any intermediate model structure is in particular fibrant in the local 
projective structure and we are done.
\end{proof}

Let $\mathcal M$ be a monoidal model category with monoidal product $\wedge$ and let $TC$ 
be the class of trivial
cofibrations in $\mathcal M$.  Recall that the monoid axiom is the statement
that all maps in $(TC \wedge \mathcal M)$-cell are weak equivalences, where
$X$-cell denotes the closure under transfinite compositions of pushouts of maps
in $X$.  This axiom ensures that the categories of modules and algebras over a
monoid in $\mathcal M$ inherit model structures from $\mathcal M$; we will
elaborate somewhat on this in Section \ref{alg-mod}.  When the category
$\mathcal M$ is cofibrantly generated with generating trivial cofibrations
$J$, then to show that the monoid axiom holds it suffices to check that every 
map in $(J \wedge \mathcal M)$-cell is a
weak equivalence.  See \cite{schwede-shipley} for further details.

In general, the intermediate model structures on $\spc$ are not monoidal, but we have the following 
important examples.

\begin{prop}
\label{spc-inj-monoidal}
The local injective model structure on $\spc$ is monoidal.
\end{prop}
\begin{proof}
Note that when $X \to Y$ is a monomorphism, so is $X \wedge Z \to Y \wedge Z$.
Given two monomorphisms $X_i \to Y_i$, consider the pushout diagram
\begin{equation}
\label{spc-inj-pushout}
\xymatrix{
X_1 \wedge X_2 \ar[r] \ar[d] & Y_1 \wedge X_2 \ar[d] \\
X_1 \wedge Y_2 \ar[r] & P \,.
}
\end{equation}
By evaluating in sections and quoting the corresponding result about simplicial sets, we get that
the induced pushout product map $P \to Y_1 \wedge Y_2$ is a monomorphism.

If in addition $X_1 \to Y_1$ is a local weak equivalence, so is the left vertical map in the pushout diagram, by
\ref{spc-smash}.  Left properness of the local injective model 
structure implies that the bottom map in \ref{spc-inj-pushout} is a local weak equivalence, and
using the 2-out-of-3 axiom we conclude that $P \to Y_1 \wedge Y_2$ is a local weak equivalence, so $\spc$ is monoidal.
\end{proof}

\begin{prop}
\label{spc-proj-monoid}
Let $\site$ be the Nisnevich site of smooth schemes of finite type over a
finite-dimensional base scheme $S$, and consider the category of spaces with the
local projective model structure.  In this case $\spc$ is a monoidal model category.
\end{prop}
\begin{proof}
First we note that the projective cofibrations are part of a monoidal model structure on spaces where the fibrations and weak equivalences are defined sectionwise, by \cite[2.7]{dro}.  To show that the local projective model structure
is monoidal it suffices to prove that the pushout product of a projective cofibration $i$ and a projective locally trivial cofibration $j$ is a local weak equivalence.  We can assume $i$ is of the form
$$
\Lambda^n_k \wedge U \to \Delta^n \wedge U
$$
and that $j \colon X \to Y$ has $X$ and $Y$ cofibrant, by the description of the generating sets in 
\cite[4.1]{blander-spre}.

In the diagram
$$
\xymatrix{
\Lambda^n_k \wedge U \wedge X \ar[d]_{1 \wedge j} \ar[r]^{i \wedge 1} & \Delta^n \wedge U \wedge X \ar[d]_{1
\wedge j} \\
\Lambda^n_k \wedge U \wedge Y \ar[r] & \Delta^n \wedge U \wedge Y
}
$$
the map $i \wedge 1$ is a cofibration since $\Lambda^n_k \to \Delta^n$ is a cofibration of simplicial sets and
$U \wedge X$ is cofibrant.  The maps denoted $1 \wedge j$ are local weak equivalences by \ref{spc-smash}.  By left
properness the induced map $\Delta^n \wedge U \wedge X \to P$, where $P$ denotes the pushout of $i \wedge 1$ and
the left $1 \wedge j$, is a local weak equivalence.  Now the pushout product $P \to \Delta^n \wedge U \wedge Y$ 
is a local weak equivalence by the 2-out-of-3 property of weak equivalences.
\end{proof}

\begin{prop}
Assume that $\spc$ is given any of the model structures constructed in Proposition \ref{spc-inter-model},
and assume in addition that the model structure is monoidal.  Then it also satisfies the monoid axiom.
\end{prop}
\begin{proof}
Consider the class $C$ of morphisms consisiting of $X \wedge Z \to Y \wedge Z$, where $Z$ is a space and 
$X \to Y$ is a trivial cofibration of spaces.  By \ref{spc-smash} this class is contained in the class of
local injective trivial cofibrations, i.e., monomorphisms which are also local weak equivalences.  Now $C$-cell
is also contained in the class of local injective cofibrations, since trivial cofibrations are closed under the 
formation of cell objects, and in particular every morphism in $C$-cell is a local weak equivalence.
\end{proof}

In Section \ref{section-stable} we will apply Bousfield localization to the
category of $\Gamma$-spaces.  For this we need to know that our categories are
combinatorial, in the sense of Jeff Smith.  An account of this notion is given
in Dugger \cite{dugger}; we recall the relevant definitions below.

\begin{defn}
Let $\lambda$ be a regular cardinal and $\cat$ a category.  An object $X \in
\cat$ is \emph{$\lambda$-presentable} if the represented functor $\cat(X,-)$
preservers $\lambda$-filtered colimits.  The category $\cat$ is \emph{locally
$\lambda$-presentable} if it is cocomplete, and there exists a set $\{G_i\}$ of
$\lambda$-presentable objects in $\cat$ such that every object in $\cat$ can
be written as a $\lambda$-filtered colimit of the $G_i$'s.  $\cat$ is
\emph{locally presentable} if it is locally $\lambda$-presentable for some
$\lambda$.
\end{defn}

\begin{defn}
A model category is \emph{combinatorial} if it is locally presentable and
cofibrantly generated.
\end{defn}

\begin{rem}
There is another notion which assures the applicability of Bousfield localization
developed in Hirschhorn's book \cite{hirschhorn}, called cellularity, which is more suitable for
categories built from topological spaces.
\end{rem}

The following basic result is found in e.g. \cite[1.12]{ar}.

\begin{prop}
\label{loc-pres}
Let $\mathcal I$ be a small category.  Then the functor category $\fun(\mathcal I,\set)$
is locally presentable.
\end{prop}

Since $\spc$ is isomorphic to $\fun(\site^{\op} \times \Delta^{\op},\set)$, we have the following result.

\begin{cor}
\label{prop-combinatorial}
The category of spaces is combinatorial.
\end{cor}

%
%
\section{The category of $\Gamma$-spaces}

Let $\Gamma$ be the full subcategory of the category of pointed sets with
objects $n_+ = \{0,1,\dots,n\}$, for $n\geq 0$,  where $0$ is the basepoint in
$n_+$.  Let $\cat$ be a pointed category.  The full subcategory of
$\fun(\Gamma,\cat)$ consisting of functors that send $0_+$ to the basepoint in $\cat$ is
the category of $\Gamma$-objects in $\cat$, denoted $\Gamma\cat$.  
Since $\Gamma$ is a skeleton for $\fset$, the category of finite based sets, we could also define
$\Gamma$-objects in $\cat$ to be the full subcategory of $\fun(\fset,\cat)$ consisting of pointed
functors.

When $\cat$
is the category $\sset$ of pointed simplicial sets, objects in $\gsset$ are classically
called $\Gamma$-spaces; model structures on this category are constructed in
the papers Bousfield and Friedlander \cite{bousfield-friedlander} and Schwede \cite{schwede}.  
Our objects of study will be $\Gamma$-objects in $\spc$, which we also call $\Gamma$-spaces.
Alternatively, our $\Gamma$-spaces can be thought of as presheaves of ordinary
$\Gamma$-spaces, i.e., $\fun(\site^{op},\gsset)$.  Note that when $\site$
consists of one morphism only, we recover the category $\gsset$, and our stable
model structure will be constructed so that we recover the stable model structure
in \cite{schwede}. 

To start with, we want to define a closed symmetric monoidal structure on $\gspc$.  Observe
that $\Gamma$ is symmetric monoidal under the operation $\wedge\colon \Gamma
\times \Gamma \to \Gamma$ given by $(m_+,n_+) \mapsto mn_+$.  
Given two $\Gamma$-spaces $F$ and $G$, the smash product $F \wedge G$
is defined as the left Kan extension filling out the following diagram.
$$ 
\xymatrix{ \Gamma \times \Gamma \ar[d]_\wedge \ar[r]^-{(F,G)} & \spc \times \spc \ar[r]^-{- \wedge -} & \spc \\ 
\Gamma \ar@{-->}[urr] }
$$ 
More explicitly, the smash product is the pointwise colimit
$$
(F \wedge G)(n_+) = \colim_{i_+ \wedge j_+ \to n_+} F(i_+) \wedge G(j_+) \,.
$$
It follows from the universal property of the colimit that maps of
$\Gamma$-spaces $F \wedge G \to H$ are in 1-1 correspondence with maps
$F(i_+) \wedge G(j_+) \to H(i_+ \wedge j_+)$ that are natural in $i_+$ and
$j_+$, and that this property characterizes $F \wedge G$ up to isomorphism.

Simplicial function complexes of $\Gamma$-spaces are defined to be
$$
\shom(F,G)_n = \gspc(F \wedge \Delta^n_+, G)
$$
in simplicial degree $n$; the face and degeneracy maps are the obvious ones.
From this we define the simplicial presheaf-hom, or space-hom, in sections by
$$
\phom(F,G)(U) = \shom(F|U,G|U) \,,
$$
where $|U$ denotes pointwise restriction to the local site $\site \downarrow U$.
Finally, internal hom-$\Gamma$-spaces are defined by setting
$$
\ihom(F,G)(n_+) = \phom(F,G(n_+ \wedge -)) \,.
$$

We have given the constructions of the objects involved in the following
result, which is a special case of Day's work in \cite{day}. 

\begin{prop}
The category $\gspc$ is a simplicial closed symmetric monoidal category
enriched over $\spc$.
\end{prop}

A set defines a discrete simplicial set, and therefore a constant simplicial
presheaf.  In particular, the sets $\Gamma(n_+,k_+)$ define the corepresented
$\Gamma$-space $\Gamma^n$ given pointwise by $\Gamma^n(k_+) = \Gamma(n_+,k_+)$.
Let $F$ be a $\Gamma$-space and let $F \circ \Gamma^n$ denote the $\Gamma$-space given pointwise by 
$$
(F \circ \Gamma^n)(k_+) = F(\Gamma(n_+,k_+)) \,.
$$

The following two lemmas follow immediately from Lydakis' corresponding results for classical $\Gamma$-spaces
in \cite{lydakis} by evaluating in sections.

\begin{lem}
\label{lem-yoneda}
There are natural isomorphisms
\begin{enumerate}
\item $\phom(\Gamma^n,F) \cong F(n_+)$
\item $\Gamma^m \wedge \Gamma^n \cong \Gamma^{mn}$
\item $F \wedge \Gamma^n \cong F \circ \Gamma^n$ \,.
\end{enumerate}
\end{lem}

\begin{lem}
\label{lem-mono}
Smashing with a $\Gamma$-space preserves monomorphisms of $\Gamma$-spaces.
\end{lem}

There are functors
\begin{equation}
\label{ev-lev}
\lev : \spc \rightleftarrows \gspc :\ev
\end{equation}
for each $n\geq 0$, where $\ev$ is evaluation at $n_+$ and $\lev(X) = X \wedge
\Gamma^n$.  From Lemma \ref{lem-yoneda} we have a natural isomorphism
\begin{equation}
\label{lev-iso}
\operatorname{L}_m(X) \wedge \lev(Y) \cong \levmn(X \wedge Y) \,.
\end{equation}

\begin{prop}
The functors in \ref{ev-lev} form an adjoint pair.
\end{prop}
\begin{proof}
Since smash products with spaces and colimits are computed pointwise, and
$\Gamma^n$ is discrete, there is a natural isomorphism 
$$
\lev(X)(k_+) \cong \bigvee_I X \,,
$$
where the index set $I$ consists of all maps $n_+ \to k_+$ in $\Gamma$ except 
for the zero map.  The natural bijection
$$
\gspc(\lev(X),G) \to \spc(X,\ev(G))
$$
is given by restricting $\lev(X)(n_+)$ to the wedge summand indexed by the
identity map $n_+ \to n_+$.  The inverse map takes $X \to G(n_+)$ to the
functor $\lev(X) \to G$ which pointwise at $k_+$ is determined on each wedge
summand $X$, indexed by $n_+ \to k_+$, by composing $X \to G(n_+)$ with $G(n_+)
\to G(k_+)$.
\end{proof}

%
%
\section{Pointwise model structures}

In this section we establish basic results about the pointwise projective model
structures on $\gspc$.  

\begin{hyp}
\label{spc-hyp}
For the rest of this paper we will assume, unless otherwise noted, that $\spc$ is given one of 
the intermediate model structures described in Theorem \ref{spc-inter-model}, including the 
local injective and local projective structures.  Let $I$ and $J$ denote the sets of generating cofibrations
and generating trivial cofibrations, respectively.  We will further assume that the domains of the maps in
$I$ and $J$ are cofibrant.
\end{hyp}

\begin{defn}
A map $F \to G$ of $\Gamma$-spaces is a 
\begin{itemize}
\item \emph{pointwise weak equivalence} if $F(n_+) \to G(n_+)$ is a local weak
equivalence in $\spc$ for all $n\geq 0$.
\item \emph{pointwise fibration} if $F(n_+) \to G(n_+)$ is a fibration in
$\spc$ for all $n\geq 0$.
\item \emph{cofibration} if it has the left lifting property with respect to
the maps that are both pointwise weak equivalences and projective fibrations.
\end{itemize}
\end{defn}

\begin{thm}
\label{thm-proj-model}
Let $I$ and $J$ be the sets of generating cofibrations and generating trivial
cofibrations in $\spc$.  Then $\gspc$ with the classes of pointwise weak
equivalences, cofibrations and pointwise fibrations is a cofibrantly generated
proper $\spc$-model category, with generating cofibrations
$$
I_\Gamma = \bigcup_{n \geq 0} \lev(I)
$$
and generating trivial cofibrations
$$
J_\Gamma = \bigcup_{n \geq 0} \lev(J) \,.
$$
We will refer to this model structure as the \emph{pointwise model structure}
on $\gspc$.
\end{thm}
\begin{proof}
This result is an application of more general results concerning pointwise
projective model structures on diagram categories, which can be found in
Hirschhorn's book, \cite[11.6.1, 11.7.3, 13.1.14]{hirschhorn}.  The model structure is
enriched in $\spc$ by \cite[3.30]{barwick}.
\end{proof}

\begin{cor}
\label{cor-proj-inj}
The adjoint functor pair \ref{ev-lev} is a Quillen pair, and $\ev$ preserves
cofibrations.  In particular, cofibrations are monomorphisms.
\end{cor}
\begin{proof}
The first statement follows immediately from Theorem \ref{thm-proj-model}, 
the second statement follows from \cite[11.6.3]{hirschhorn}.
\end{proof}

\begin{cor}
\label{cofibrant-gspaces}
The $\Gamma$-space $X \wedge \Gamma^n$ is cofibrant when $X$ is a cofibrant
space.  In particular $\Gamma^n$ is cofibrant.
\end{cor}
\begin{proof}
This follows by applying $\lev$ to the map $* \to X$.
\end{proof}

As $\gspc$ as a category is isomorphic to 
$\fun(\Gamma \times \site^{\op} \times \Delta^{\op}, \set)$, it is
locally presentable by Proposition \ref{loc-pres}.

\begin{cor}
The category of\/ $\Gamma$-spaces with the pointwise model structure is
combinatorial.
\end{cor}

\begin{prop}
\label{prop-proj-monoidal}
The category of\/ $\Gamma$-spaces equipped with the pointwise model structure
is a monoidal model category provided $\spc$ is monoidal.
\end{prop}
\begin{proof}
Since the monoidal unit $\Gamma^1$ is cofibrant, it suffices to check the pushout product
axiom.  Let $F_i \to G_i$, where $i=1,2$, be two cofibrations, and construct the pushout diagram
\begin{equation}
\label{pushout-product}
\xymatrix{
F_1 \wedge F_2 \ar[r] \ar[d] & G_1 \wedge F_2 \ar[d] \\
F_1 \wedge G_2 \ar[r] & P \,.
}
\end{equation}
We need to show that the induced pushout product map $P \to G_1 \wedge G_2$ is
a cofibration.  We may assume the $F_i \to G_i$ are of the form 
$$
X_i \wedge \Gamma^{n_i} \to Y_i \wedge \Gamma^{n_i} \,,
$$
where $X_i \to Y_i$ are cofibrations in $\spc$.  Using the isomorphism
\ref{lev-iso}, and the fact that $\levnn$ preserves colimits, we can apply
$\levnn$ to the pushout constructed from the maps $X_i \to Y_i$ to obtain
$$
\xymatrix{
\levnn(X_1 \wedge X_2) \ar[r] \ar[d] & \levnn(Y_1 \wedge X_2) \ar[d] \\
\levnn(X_1 \wedge Y_2) \ar[r] & \levnn(X_1 \wedge Y_2 \coprod_{X_1 \wedge X_2} Y_1 \wedge X_2) \,,
}
$$
which is isomorphic to \ref{pushout-product}.  We know that
$$
X_1 \wedge Y_2 \coprod_{X_1 \wedge X_2} Y_1 \wedge X_2 \to X_2 \wedge Y_2
$$
is a cofibration of spaces, by the assumption that $\spc$ is monoidal, so
the map $P \to \levnn(Y_1 \wedge Y_2)$ is a
cofibration.  The same argument gives the corresponding result about trivial
cofibrations.
\end{proof}

\begin{prop}
The pointwise model structure on $\gspc$ satisfies the monoid axiom when $\spc$
does.
\end{prop}
\begin{proof}
We need to show that the maps in $(J_\Gamma \wedge \gspc)$-cell are weak
equivalences. Consider first a map $f$ of the form
$$
\lev(X) \wedge F \to \lev(Y) \wedge F
$$
where $X \to Y$ is a generating trivial cofibration in $\spc$.  Evaluating at
$k_+$, we get
$$
\xymatrix{ X \wedge (\Gamma^n \wedge F)(k_+) \ar[r]^{f(k_+)} & Y \wedge
(\Gamma^n \wedge F)(k_+) \,, }
$$
so $f(k_+)$ is in $J \wedge \spc$ for all $k_+$.  Now, if $g$ is in $(J_\Gamma
\wedge \gspc)$-cell, it is a transfinite composition of pushouts of maps $f_i$
in $J_\Gamma \wedge \gspc$.  Since each $f_i(k_+)$ is in $J \wedge \spc$, and
colimits in $\gspc$ are computed pointwise, $g(k_+)$ is in $(J \wedge
\spc)$-cell.  Using the assumption that the monoid axiom holds in $\spc$, we
see that $g(k_+)$ is a weak equivalence for all $k_+$.
\end{proof}

\begin{lem}
\label{lem-colim}
A filtered colimit of pointwise equivalences is a pointwise equivalence.
\end{lem}
\begin{proof}
Given a filtered category $\mathcal I$, consider the colimit functor
$$
\colim \colon \fun({\mathcal I},\gspc) \to \gspc \,.
$$
This is a left Quillen functor, where $\fun({\mathcal I},\gspc)$ is given the pointwise projective model structure
(see the proof of Theorem \ref{thm-proj-model}.)  A weak equivalence in $\fun(\mathcal I,\gspc)$ factors as a trivial
cofibration followed by a trivial fibration, and since $\colim$ preserves trivial cofibrations, it remains to show
that $\colim$ preserves trivial fibrations.

Let $\{F_\alpha\} \to \{G_\alpha\}$ be a trivial fibration in $\fun(\mathcal I,\gspc)$, and consider the
lifting problem
$$
\xymatrix{
A \ar[r] \ar[d] & \colim F_\alpha \ar[d] \\
B \ar[r] & \colim G_\alpha
}
$$
where $A \to B$ is a generating cofibration in $\gspc$.  The existence of a lift in this diagram is equivalent to the surjectivity of 
the canonical map
\begin{equation}
\label{lift-map}
\gspc(B,\colim F_\alpha) \to \gspc(A,\colim F_\alpha) \times_{\gspc(A,\colim G_\alpha)} \gspc(B,\colim G_\alpha) \,.
\end{equation}
Since $A$ and $B$ are small and finite limits commute with filtered colimits in the category of sets, 
\ref{lift-map} is the colimit of the canonical maps
\begin{equation}
\label{lift-map2}
\gspc(B,F_\alpha) \to \gspc(A,F_\alpha) \times_{\gspc(A,G_\alpha)} \gspc(B,G_\alpha) \,.
\end{equation}
All the maps in \ref{lift-map2} are surjective, since this corresponds to lifts in diagrams of $\Gamma$-spaces involving cofibrations
and trivial fibrations.  A filtered colimit of surjective maps of sets is surjective, and we are done.
\end{proof}

\begin{rem}
More generally one can consider the class of ``weakly finitely generated'' model categories in the sense of
Dundas, R\"ondigs and {\O}stv{\ae}r \cite[3.4]{dro-en}.  These model categories in particular have the property 
stated in Lemma \ref{lem-colim}, and our proof of \ref{lem-colim} is taken from \cite[3.5]{dro-en}.
\end{rem}

\begin{prop}
\label{pt-smash}
Pointwise equivalences of $\Gamma$-spaces are preserved when smashed with a cofibrant $\Gamma$-space.
\end{prop}
\begin{proof}
Let $f \colon F \to G$ be a pointwise equivalence.  The induced map $F \circ \Gamma^n \to G \circ \Gamma^n$
is clearly a pointwise equivalence, so by Lemma \ref{lem-yoneda}, the map 
$f \wedge 1 \colon F \wedge \Gamma^n \to G \wedge \Gamma^n$ is a pointwise equivalence.  

Now let $C$ be a cofibrant $\Gamma$-space.  Since $\gspc$ is cofibrantly
generated with generating cofibrations $I_\Gamma$, $C$ is a retract of an $I_\Gamma$-cell complex, where by 
$I_\Gamma$-cell complex we mean that the unique map $* \to C$ is a transfinite composition of pushouts of maps in $I_\Gamma$.  Weak equivalences are closed under retracts, so it suffices to consider $C = \colim_{\alpha<\gamma} C_\alpha$, $\gamma$ an ordinal, where the maps $C_\alpha \to C_{\alpha+1}$ are given by pushout diagrams
\begin{equation}
\label{cell-pushout}
\xymatrix{
X \wedge \Gamma^n \ar[r] \ar[d]_{i \wedge 1} & C_\alpha \ar[d] \\
Y \wedge \Gamma^n \ar[r] & C_{\alpha+1} \,.
}
\end{equation}
Here $i \colon X \to Y$ is a cofibration of spaces.

Smashing \ref{cell-pushout} with $F$ and $G$ gives us two pushout diagrams as the top and bottom faces of a cubical
diagram.  Assuming by induction that $F \wedge C_\alpha \to G \wedge C_\alpha$ is a pointwise equivalence, the
gluing lemma (see \cite[II.8.12]{gj}) can be applied to conclude that $F \wedge C_{\alpha+1} \to G \wedge C_{\alpha+1}$ is a 
pointwise equivalence.  Since $F \wedge C \to G \wedge C$ is the colimit of the maps $F \wedge C_\alpha \to G \wedge C_\alpha$ 
we can conclude by applying Lemma \ref{lem-colim}.
\end{proof}

%
%
\section{Stable model structures}
\label{section-stable}

In this section we will construct the stable model structures for (presheaves of) $\Gamma$-spaces and compare it
to the model category of (presheaves of) spectra.  In fact, parts of our construction relies on this comparison; we
will begin by recalling the theory of spectra on a site.

For us, a spectrum is a sequence of objects $E^k \in \spc$ indexed by non-negative integers $k$ together
with structure maps
$$
S^1 \wedge E^k \to E^{k+1}
$$
for each $k$.  Maps of spectra are sequences of maps $f^k\colon E^k \to F^k$
compatible with the structure maps in the sense that the diagram
$$
\xymatrix{
S^1 \wedge E^k \ar[r] \ar[d]_{1 \wedge f^k} & E^{k+1} \ar[d]^{f^{k+1}} \\
S^1 \wedge F^k \ar[r] & F^{k+1}
}
$$
commutes for all $k$.  Denote the category of spectra by $\spt$.

A spectrum $E$ is levelwise fibrant if each $E^k$ is fibrant, and is
an $\Omega$-spectrum if the adjoints $E^k \to \Omega E^{k+1}$ of the structure maps
are weak equivalences.  The loop functor $\Omega\colon \spc \to \spc$ is by definition a fibrant replacement $(-)_f$ followed by the simplicial cotensor $(-)^{S^1}$ on spaces.  Note that we do not require our $\Omega$-spectra to be
levelwise fibrant.  A map $f \colon E \to F$ of spectra is a cofibration if $f^0 \colon E^0 \to F^0$ is a cofibration
of spaces and the induced maps 
$$
(S^1 \wedge F^k) \bigcup_{S^1 \wedge E^k} E^{k+1} \to F^{k+1}
$$
are cofibrations of spaces for all $k\geq 0$.
The map $f$ is a stable equivalence of spectra if it induces isomorphisms $\pi_n(E) \to \pi_n(F)$ of stable homotopy sheaves
for all integers $n$ and $U \in \site$, where the stable homotopy sheaf $\pi_n(E)$ is by definition the colimit of
the system
$$ 
\dots \to \pi_{n+k}(E^k) \to \pi_{n+k+1}(S^1 \wedge E^k) \to \pi_{n+k+1}(E^{k+1}) \to \dots \,.
$$

The following result was first proved by Jardine in \cite[2.8]{jardine-spt} for the local injective model structure on 
$\spc$; Hovey has results for spectra in more general model categories in \cite[3.3]{hovey}.

\begin{thm}
With the above notions of stable cofibrations and stable equivalences the category $\spt$ of
spectra is a cofibrantly generated proper $\spc$-model category.
A spectrum is stably fibrant if and only if it is a levelwise fibrant $\Omega$-spectrum.
\end{thm}

Let $F$ be a $\Gamma$-space, which we now consider as a based functor from all finite based sets to $\spc$.
The functor $F$ induces a functor $\wF \colon \sset \to \sspc$ from simplicial
sets to simplicial spaces, by applying $F$ in each simplicial degree.  We can
compose $F$ with the diagonal functor $d \colon \sspc \to \spc$ to get a functor 
$$
d\wF \colon \sset \to \spc \,.
$$

\begin{prop}
\label{F-eq}
Let $K \to L$ be a weak equivalence of simplicial sets.  Then the induced map 
$d\wF(K) \to d\wF(L)$ is a sectionwise equivalence, and in particular a local weak equivalence.
\end{prop}
\begin{proof}
This follows from the corresponding result for classical $\Gamma$-spaces in \cite[4.9]{bousfield-friedlander}, since
$d\wF(K)(U)$ coincides with the corresponding construction for the classical $\Gamma$-space $F(U)$.
\end{proof}

Each pair of based sets $U,V$ induces natural maps
$$
U \wedge F(V) \to F(U \wedge V)
$$
whose adjoints $U \to \spc(F(V),F(U \wedge V))$ are described by sending an element
$u$ to the map $F(u \wedge -)$.  These maps induce simplicial maps
$$
X \wedge \wF(Y) \to \wF(X \wedge Y)
$$
where $X$ and $Y$ are based simplicial sets.  By applying the diagonal functor this results in
maps
\begin{equation}
\label{gspc-sp-maps}
X \wedge d\wF(Y) \to d\wF(X \wedge Y) \,.
\end{equation}

The spectrum associated to a $\Gamma$-space $F$, which we denote $\Sp(F)$, is
defined on each level as $\Sp(F)^n = d\wF(S^n)$.  Here $S^n = S^1 \wedge \dots \wedge S^1$ ($n$ times.)
As a special case of \ref{gspc-sp-maps} we have
$$
S^m \wedge d\wF(S^n) \to d\wF(S^{m+n})
$$
which gives us the structure maps for $\Sp(F)$.

\begin{lem}
\label{sp-props}
The functor $\Sp(F)$ has the following properties.
\begin{enumerate}
\item $\Sp(F)^0 = F(1_+)$
\item $\Sp(\Gamma^n) = \s^{\times n}$
\item $\Sp(X \wedge F) = X \wedge \Sp(F)$, for spaces $X$.
\end{enumerate}
\end{lem}

Let $E$ be a spectrum.  We obtain a $\Gamma$-space $\gae$ by defining
$$
\gae(n_+) = \phom_{\spt}(\s^{\times n}, E) \,,
$$
where $\s$ denotes the sphere spectrum.  Here $\phom_{\spt}(-,-)$ denotes the
space of morphisms in the category of spectra, defined sectionwise in the same way as for
$\Gamma$-spaces, i.e.,
$$
\phom_{\spt}(E,F)(U) = \shom_{\spt}(E|U,F|U)
$$
for all $U \in \site$.  A morphism $\theta \colon m_+ \to n_+$ induces a map
$\theta^* \colon \s^{\times n} \to \s^{\times m}$ by copying the $\theta(i)$'th
factor into the $i$'th factor.  This map in turn induces $\gae(m_+) \to \gae(n_+)$.

\begin{lem}
\label{S-coeq}
The spectrum $\Sp(F)$ coincides with the coequalizer of the diagram
$$
\xymatrix{
\underset{\theta\colon m_+ \to n_+} \bigvee \s^{\times n} \wedge F(m_+) 
\ar@<1ex>[r]^-{1 \wedge F(\theta)} \ar[r]_-{\theta^* \wedge 1} & 
\underset{k_+}\bigvee \s^{\times k} \wedge F(k_+) \,.
}
$$
\end{lem}
\begin{proof}
Since colimits in $\spt$, $\spc$ and $\sset$ are computed pointwise, 
it suffices to show that the following diagram
$$
\xymatrix{
\underset{\theta\colon m_+ \to n_+}\bigvee (S^i_q)^{\times n} \wedge F(m_+)(U)_q
\ar@<1ex>[r]^-{1 \wedge F(\theta)} \ar[r]_-{\theta^* \wedge 1} & 
\underset{k_+}\bigvee (S^i_q)^{\times k} \wedge F(k_+)(U)_q \ar[r]^-f & F(S^i_q)(U)_q \,.
}
$$
is a coequalizer of sets, for all $i,q\geq 0$, where $f$ is
described as follows.  A collection of $k$ ordered elements $x_j$ in $S^i_q$
specifies a map $k_+ \to S^i_q$, and by applying $F$ we get a map 
$(S^i_q)^{\times k} \to \set(F(k_+)(U)_q, F(S^i_q)(U)_q)$.  Take the adjoint of this
and sum over $k_+$ to get $f$.  We omit the straightforward element chase.
\end{proof}

\begin{prop}
\label{prop-adj}
The functors
$$
\Sp : \gspc \rightleftarrows \spt : \Phi
$$
constitute an adjoint pair.  Furthermore, this adjunction can be extended to a $\spc$-adjunction
$$
\phom_{\spt}(\Sp(F),E) \cong \phom_{\gspc}(F,\Phi(E)) \,.
$$
\end{prop}
\begin{proof}
First note that we have an adjunction
\begin{equation}
\label{yet-another-adj}
\spt(X \wedge E,F) \cong \spc(X,\phom(E,F))
\end{equation}
where $X$ is a space and $E,F$ are spectra.  Now take the coequalizer in Lemma \ref{S-coeq}
and apply the functor $\spt(-,E)$ and the isomorphism \ref{yet-another-adj}.  The result
is that $\spt(\Sp(F),E)$ is the equalizer of
$$
\xymatrix{
\underset{\theta\colon m_+ \to n_+}\prod \spc(F(m_+),\phom(\s^{\times n},E)) & 
\underset{k_+}\prod \spc(F(k_+),\phom(\s^{\times k},E)) \ar@<-1ex>[l] \ar[l] \,.
}
$$
Any map $\Sp(F) \to E$ thus corresponds to a collection of maps 
$$
F(k_+) \to \phom(\s^{\times k},E) = \Phi(E)(n_+)
$$
natural in $k_+$, i.e., a map of $\Gamma$-spaces $F \to \Phi(E)$.
\end{proof}

\begin{prop}
\label{s-preserve}
The functor $\Sp$ preserves cofibrations.
\end{prop}
\begin{proof}
It suffices to consider generating cofibrations in $\gspc$.
Let $X \wedge \Gamma^n \to Y \wedge \Gamma^n$ be a generating cofibration, where $X \to Y$ is a cofibration of spaces.
By Lemma \ref{sp-props} we need to show that
$$
X \wedge \Sp(\Gamma^n) \to Y \wedge \Sp(\Gamma^n)
$$
is a cofibration of spectra, but this is immediate since $\Sp(\Gamma^n) = \s^{\times n}$ is a cofibrant spectrum and
$\spt$ is a $\spc$-model category.
\end{proof}

\begin{defn}
\label{defn-special}
A $\Gamma$-space $F$ is \emph{special} if the maps
$$
F(n_+) \to F(1_+) \times \dots \times F(1_+)
$$
induced by the $n$ usual projections from $n_+$ to $1_+$ are weak equivalences for
all $n \geq 1$.  If, in addition, the map
$$
F(2_+) \to F(1_+) \times F(1_+)
$$
induced by a projection and the fold map is a weak equivalence, then $F$ is
\emph{very special}.
\end{defn}

Note that when $F$ is special the maps
$$
\xymatrix{F(1_+) \times F(1_+) & F(2_+) \ar[l]_-\sim \ar[r]^-\nabla & F(1_+)}
$$
induce a commutative monoid structure on $\pi_0(F(1_+))$.  If $F$ is very special, then $\pi_0(F(1_+))$ is in fact 
an abelian group.

\begin{prop}
\label{s-phi-preserve}
The functor $\Sp$ sends very special $\Gamma$-spaces to $\Omega$-spectra.
The functor $\Phi$ sends fibrant spectra to pointwise fibrant very special $\Gamma$-spaces.
\end{prop}
\begin{proof}
Let $F$ be a very special $\Gamma$-space, and let $F \to F_f$ be a pointwise fibrant replacement.
Since $\wF(S^n) \to \wF_f(S^n)$ is a pointwise equivalence of simplicial spaces, the induced map
$d\wF(S^n) \to d\wF_f(S^n)$ of spaces is a local equivalence by Proposition \ref{diagonal-eq}.  We need to show that
the map $d\wF_f(S^n) \to \Omega d\wF_f(S^{n+1})$ is a local weak equivalence.  Since $F$ is very special, so is $F_f$,
and in fact the maps 
$$
F_f(n_+) \to F_f(1_+) \times \dots \times F_f(1_+)
$$ 
and 
$$
F_f(2_+) \to F_f(1_+) \times F_f(1_+)
$$
are sectionwise equivalences by Proposition \ref{spc-loc} since each $F_f(n_+)$ is fibrant.  Thus $F_f(U)$ is a very 
special 
$\Gamma$-space in the classical sense, for each $U \in \site$, and by \cite[4.2]{bousfield-friedlander} each map
$d\wF_f(U)(S^n) \to \Omega d\wF_f(U)(S^{n+1})$ is a weak equivalence of simplicial sets.  This implies in particular
that $d\wF_f(S^n) \to \Omega d\wF_f(S^{n+1})$ is a local weak equivalence.

For the second statement, let $E$ be a fibrant spectrum.  Since $\s \vee \dots \vee \s \to \s \times \dots \times \s$ is a stable equivalence of cofibrant spectra, the map
$$
\phom(\s^{\times n},E) \to \phom(\s^{\vee n},E) \cong \phom(\s,E)^{\times n}
$$
is a local weak equivalence, i.e., $\Phi(E)$ is special.  Similarly, the map $\s \vee \s \to \s \times \s$ induced by
an inclusion and the diagonal map is a stable equivalence, so $\Phi(E)$ is very special.
\end{proof}

\begin{defn}
The \emph{$n$-th homotopy sheaf} $\pi_n(F)$ of a $\Gamma$-space $F$ is the
$n$-th homotopy sheaf of the associated spectrum $\Sp(F)$.
We write $\pi_*(F)$ for the $\Z$-graded abelian sheaf $\oplus_n \, \pi_n(F)$.
\end{defn}

Note that an equivalent definition of $\pi_n(F)$ is as the sheaf associated to the presheaf 
$U \mapsto \pi_n(F(U))$, where $\pi_n(F(U))$ are homotopy groups of classical $\Gamma$-spaces.

\begin{defn}
\label{defn-stable}
A map $F \to G$ in $\gspc$ is a
\begin{itemize}
\item \emph{stable equivalence} if the induced map $\pi_*(F) \to \pi_*(G)$ is an isomorphism.
\item \emph{stable fibration} if it has the right lifting property with respect to the maps that are both cofibrations and stable equivalences.
\end{itemize}
\end{defn}

Recall that a spectrum $E$ is called connective if $\pi_n(E)=0$ for $n<0$.
Since the $k$-simplices of $\Delta^n/\partial\Delta^n$ for $k<n$ consist of the basepoint only, and since 
$\Delta^n/\partial\Delta^n$ is weakly equivalent to $S^n$, it follows from Proposition \ref{F-eq} that $dF(S^n)$ is
$(n-1)$-connected,  and that $\Sp(F)$ is a connective spectrum.

\begin{lem}
\label{unit-counit}
The following holds for the adjunction in Proposition \ref{prop-adj}.
\begin{enumerate}
\item The composition $F \to \Phi(\Sp(F)) \to \Phi(\Sp(F)_f)$ of the unit map and $\Phi$ applied to a fibrant replacement of 
$\Sp(F)$, is a pointwise weak equivalence for special $\Gamma$-spaces $F$.
\item When $E$ is a fibrant spectrum, the counit map $\Sp(\Phi(E)) \to E$ induces isomorphisms 
$\pi_n(\Sp(\Phi(E))) \to \pi_n(E)$ for all $n\geq 0$.  In particular $\Sp(\Phi(E)) \to E$ is a stable equivalence when $E$ is a
fibrant connective spectrum.
\end{enumerate}
\end{lem}
\begin{proof}
Let $F$ be special.  The commutative diagram
$$
\xymatrix{
F(n_+) \ar[r] \ar[d]^\sim & \phom(\s^{\times n},\Sp(F)_f) \ar[d]^\sim \\
F(1_+)^{\times n} \ar[d]^\sim & \phom(\s^{\vee n},\Sp(F)_f) \ar[d]^\cong \\
(\Sp(F)_f^0)^{\times n} \ar[r]^-\cong & \phom(\s,\Sp(F)_f)^{\times n}
}
$$
shows that the top map is a local weak equivalence for each $n\geq 0$.

When $E$ is a fibrant spectrum, $\pi_n(E) \cong \pi_n(E^0)$ for all $n\geq 0$, so
the second statement of the lemma is reduced to the statement that
$\Sp(\Phi(E))^0 \to E^0$ is a local weak equivalence of spaces.  But this map
coincides with the canonical weak equivalence
$$
\Sp(\Phi(E))^0 = (\Phi E)(1_+) = \phom_{\spt}(\s,E) \to E^0 \,.
$$
\end{proof}

We let $\ho(\spt)_{\geq 0}$ denote the full subcategory of $\ho(\spt)$ consisting of the connective spectra.

\begin{thm}
\label{thm-stable-model}
The category $\gspc$ with the classes of stable equivalences, cofibrations and
stable fibrations is a cofibrantly generated left proper $\spc$-model category,
such that the functor pair in Proposition \ref{prop-adj} induces an equivalence of categories
$$
L\Sp : \ho(\gspc) \simeq \ho(\spt)_{\geq 0} : R\Phi \,.
$$
The stably fibrant objects in $\gspc$ are the very special $\Gamma$-spaces that are also pointwise fibrant.
A pointwise fibration of stably fibrant $\Gamma$-spaces is necessarily a stable fibration.
A stable equivalence between stable fibrant $\Gamma$-spaces is a pointwise equivalence.
\end{thm}
\begin{proof}
Let $\Sigma$ be the set of maps consisting of
$$
\Gamma^1 \vee \dots \vee \Gamma^1 \to \Gamma^n
$$
for all $n\geq 1$, and the shear map
$$
\Gamma^1 \vee \Gamma^1 \to \Gamma^2 \,.
$$
These morphisms are induced by the same morphisms in $\Gamma$ as in Definition \ref{defn-special}, and corepresent the
morphisms displayed there. 
Since the pointwise model structure on $\gspc$ is combinatorial, left proper and enriched over $\spc$, we can apply enriched left Bousfield localization (see Theorem \ref{enriched-bousfield}) with 
respect to $\Sigma$ to obtain a new combinatorial and left proper model structure on $\gspc$.  For the remainder of this
proof we will refer to this model structure as the ``localized model structure.''

The localized fibrant objects are given by the $\Sigma$-local objects.  A $\Gamma$-space $H$ is $\Sigma$-local if and 
only if it is pointwise fibrant and the maps
$$
\phom(\Gamma^n,H) \to \phom(\Gamma^1 \vee \dots \vee \Gamma^1,H)
$$
and
$$
\phom(\Gamma^2,H) \to \phom(\Gamma^1 \vee \Gamma^1,H)
$$
are weak equivalences of spaces, for $n\geq 1$.  Composing with the isomorphism
$$
\phom(\Gamma^1 \vee \dots \vee \Gamma^1,H) \to \phom(\Gamma^1,H) \times \dots \times \phom(\Gamma^1,H)
$$
and using the isomorphism (1) in Lemma \ref{lem-yoneda}, it is
clear that the $\Sigma$-local objects coincide with the pointwise fibrant very special $\Gamma$-spaces.

The localized weak equivalences are defined to be those maps $f\colon F \to G$ that have a cofibrant replacement
$f_c\colon F_c \to G_c$ (in the pointwise model structure) that induces local weak equivalences
$$
\phom(G_c,H) \to \phom(F_c,H)
$$
of spaces for all $\Sigma$-local $H$.  We have to identify the localized weak equivalences
as the stable equivalences.

Consider the following diagram
\begin{equation}
\label{map-diagram}
\xymatrix{
\phom(G_c,\Phi(E)) \ar[r]^-\cong \ar[d]_{f_c^*} & \phom(\Sp(G_c),E) \ar[d]^{\Sp(f_c)^*} \\
\phom(F_c,\Phi(E)) \ar[r]^-\cong & \phom(\Sp(F_c),E)
}
\end{equation}
where the horizontal maps come from the simplicial version of the adjunction in Proposition
\ref{prop-adj}.  Note that $\Sp(f_c)$ is a map between cofibrant objects by
\ref{s-preserve}.  Since $\spt$ is a simplicial model category, $\Sp(f_c)\colon
\Sp(F_c) \to \Sp(G_c)$ is a stable equivalence of spectra if and only if
$\Sp(f_c)^*$ is a weak equivalence of simplicial sets for all fibrant spectra
$E$.  It follows that $f_c$ is a stable equivalence of $\Gamma$-spaces if and
only if $f_c^*$ is a weak equivalence for all fibrant $E$.  In particular, 
a localized weak equivalence is a stable equivalence since by Proposition \ref{s-phi-preserve} 
we know that $\Phi(E)$ is a $\Sigma$-local $\Gamma$-space.

When $H$ is a very special $\Gamma$-space the map $H \to \Phi(\Sp(H)_f)$ is a pointwise weak
equivalence by \ref{unit-counit}, and hence induces weak equivalences of
simplicial sets in the diagram
\begin{equation}
\label{banankake}
\xymatrix{
\phom(G_c,H) \ar[r]^-\sim \ar[d]_{f_c^*} & \phom(G_c,\Phi(\Sp(H)_f)) \ar[d]^{f_c^*} \\
\phom(F_c,H) \ar[r]^-\sim & \phom(F_c,\Phi(\Sp(H)_f)) \,.
}
\end{equation}
It follows from \ref{map-diagram} and \ref{banankake} that a stable equivalence is a
localized weak equivalence.

Now that we have identified the localized weak equivalences as the stable
equivalences, $\Sp$ becomes a left Quillen functor by \ref{s-preserve} since the
localization process does not change the class of cofibrations.  The Quillen
pair $\Sp$ and $\Phi$ induces derived adjoint functors $L\Sp$ and $R\Phi$ on the
homotopy categories of $\gspc$ and $\spt$, which by \ref{s-phi-preserve} restrict to functors
$$
L\Sp : \ho(\gspc) \rightleftarrows \ho(\spt)_{\geq 0} : R\Phi \,.
$$
To show that $L\Sp$ is an equivalence, it is enough to note that $\Sp$ detects
weak equivalences, and that the counit map $\Sp(\Phi(E)) \to E$ is a stable
equivalence for connective fibrant spectra $E$ by Lemma \ref{unit-counit}. 
\end{proof}

\begin{prop}
\label{cobase-smash}
Smashing with a cofibrant $\Gamma$-space preserves stable equivalences.
\end{prop}
\begin{proof}
First note that $\ihom(C,H)$ is very special when $C$ is cofibrant and $H$ is
fibrant, since $\gspc$ is a $\spc$-model category. 
Let $f\colon F \to G$ be stable equivalence with cofibrant replacement $f_c\colon F_c \to G_c$, 
and $C$ a cofibrant $\Gamma$-space.  
We have that $\shom(G_c,H) \to \shom(F_c,H)$ is a
weak equivalence for all fibrant $H$, so in particular
$$
\shom(G_c,\ihom(C,H)) \to \shom(F_c,\ihom(C,H))
$$
is a weak equivalence for all cofibrant $C$ and fibrant $H$.  Together
with the isomorphism $\shom(F_c,\ihom(C,H)) \cong \shom(F_c \wedge C,H)$ this
implies that $f_c \wedge 1$ is a stable equivalence.  The commutative diagram
$$
\xymatrix{
F_c \wedge C \ar[r] \ar[d]_{f_c \wedge 1} & F \wedge C \ar[d]_{f \wedge 1} \\
G_c \wedge C \ar[r] & G \wedge C \,,
}
$$
where the horizontal maps are pointwise weak equivalences by Proposition \ref{pt-smash}, implies that 
$f \wedge 1$ is a stable equivalence.
\end{proof}

\begin{lem}
\label{lem-exactseq}
Let $F \to G$ be a monomorphism of $\Gamma$-spaces.  Then there is an exact sequence of abelian sheaves
$$
\dots \to \pi_{n+1}(G/F) \to \pi_n(F) \to \pi_n(G) \to \pi_n(G/F) \to \pi_{n-1}(F) \to \dots \,.
$$
\end{lem}
\begin{proof}
This follows from \cite[1.3]{schwede} by evaluating in sections and applying the exact sheafification functor.
\end{proof}

\begin{prop}
\label{pushout-stable}
Pushouts of $\Gamma$-spaces preserve monomorphic stable equivalences.
\end{prop}
\begin{proof}
Consider the pushout diagram
$$
\xymatrix{
F \ar[d] \ar[r] & G \ar[d] \\
F' \ar[r] & G'
}
$$
where $F \to G$ is an injective stable equivalence.  It follows that the map $F' \to G'$ is injective, and that
$G'/F' \cong G/F$, so by Lemma \ref{lem-exactseq} the map $F' \to G'$ is also a stable equivalence.
\end{proof}

\begin{prop}
\label{prop-stable-monoidal}
The stable model structure on $\gspc$ is monoidal when $\spc$ is monoidal.
\end{prop}
\begin{proof}
The first part of the pushout product axiom is immediate from Proposition
\ref{prop-proj-monoidal}.  Given a pushout diagram
$$
\xymatrix{
\lev(X_1) \wedge \levm(X_2) \ar[r] \ar[d] & \lev(Y_1) \wedge \levm(X_2) \ar[d] \\
\lev(X_1) \wedge \levm(Y_2) \ar[r] & P \,,
}
$$
it suffices to check that the induced map $P \to \lev(Y_1) \wedge \levm(Y_2)$ is
a trivial cofibration when $X_1 \to Y_1$ is a generating cofibration and
$X_2 \to Y_2$ is a generating trivial cofibration.

First note that $\lev(X_1)$ and $\lev(Y_1)$ are cofibrant.
The left vertical map in the pushout diagram is a monomorphism by \ref{lem-mono}, and a stable equivalence by
\ref{cobase-smash}.  By Proposition \ref{pushout-stable} the right vertical map is a stable equivalence; the 
pushout product
map is now seen to be a stable equivalence by the 2-out-of-3 property of stable equivalences.
\end{proof}

\begin{prop}
\label{prop-stable-monoid}
The stable model structure on $\gspc$ satisfies the monoid axiom when $\spc$ is monoidal.
\end{prop}
\begin{proof}
Let $F \to G$ be a trivial cofibration and let $H$ be a $\Gamma$-space.  The induced map
$F \wedge H \to G \wedge H$ is a monomorphism by \ref{lem-mono}, and we claim that the cofibre
$(G/F) \wedge H$ is stably contractible, which by \ref{lem-exactseq} implies that $F \wedge H \to
G \wedge H$ is a stable equivalence.  First take a cofibrant replacement $H_c \to H$.  Since $* \to G/F$ is a
stable equivalence, $(G/F) \wedge H_c$ is stably contractible by \ref{cobase-smash}, and also, $(G/F) \wedge H_c$ is
stably equivalent to $(G/F) \wedge H$, which proves the claim.

By Proposition \ref{pushout-stable}, it remains to show that a transfinite composition of stable equivalences is a 
stable equivalence.  Note first that homotopy groups of $\Gamma$-spaces commute with filtered colimits,
since this is true for spectra of simplicial sets and sheafification is exact.
A transfinite composition $F_0 \to \colim_\alpha F_\alpha$, where each $F_\alpha \to F_{\alpha+1}$ is
a stable equivalence, induces an isomorphism
$\pi_*F_0 \to \colim_\alpha \pi_*(F_\alpha) \cong \pi_*(\colim_\alpha F_\alpha)$.
\end{proof}

A symmetric spectrum is a spectrum $E$ with a $\Sigma_n$-action on each $E^n$
such that the iterated structure maps
$$
S^m \wedge E^n \to S^{m-1} \wedge E^{1+n} \to \dots \to E^{m+n}
$$
are $\Sigma_m \times \Sigma_n$-equivariant, where $\Sigma_m \times \Sigma_n$ is
identified with a subgroup of $\Sigma_{m+n}$ in the usual way.  Symmetric
spectra form a subcategory of the category of spectra, where the morphisms are 
maps of spectra equivariant at each level.  We denote this category by $\sptsym$.

Let $U \colon \sptsym \to \spt$ denote the forgetful functor, which is right
adjoint to a ``free symmetric spectrum'' functor $F \colon \spt \to \sptsym$.
A map $f \colon E \to F$ of symmetric spectra is a fibration if $U(f)\colon
U(E) \to U(F)$ is a fibration of spectra.  There are simplicial mapping spaces
of symmetric spectra, and weak equivalences of symmetric spectra are those maps
$f$ which induce weak equivalences of simplicial sets $\shom(F,H) \to
\shom(E,H)$ for all fibrant symmetric spectra $H$.  If $U(f)$ is a stable
equivalence of spectra, then $f$ is a weak equivalence of symmetric spectra,
but the converse is not true.

The following theorem is a special case of a result by Hovey \cite[8.7]{hovey}.

\begin{thm}
With the above definitions of fibrations and stable equivalences $\sptsym$ is a
cofibrantly generated proper $\spc$-model category, such that
$$
F : \spt \rightleftarrows \sptsym : U
$$
defines a Quillen equivalence.
\end{thm}

As the $\Sigma_n$-action on $S^n$ induces an action on $d\wF(S^n)$, the functor
$\Sp$ factors through the category of symmetric spectra in the sense that we have
a commutative diagram
$$
\xymatrix{
\gspc \ar[rr]^\Sp \ar[rd]_{\Sp^\Sigma} & & \spt \\
{} & \sptsym \ar[ru]_U \,.
}
$$

\begin{prop}
The functor $\Sp^\Sigma$ is lax monoidal.
\end{prop}
\begin{proof}
We can just evaluate in sections and apply the corresponding result for classical $\Gamma$-spaces and symmetric spectra in \cite[3.3]{mmss}.
\end{proof}

Note that $\Sp^\Sigma$ is not strict monoidal since 
$\Sp^\Sigma(\Gamma^m \wedge \Gamma^n) = \Sp^\Sigma(\Gamma^{mn}) = \s^{\times mn}$,
while $\Sp^\Sigma(\Gamma^m) \wedge \Sp^\Sigma(\Gamma^n) = \s^{\times m} \wedge \s^{\times n}$.
Nor is $\Sp^\Sigma$ a left Quillen functor, as $\Sp^\Sigma(\Gamma^n) = \s^{\times n}$ is not a cofibrant symmetric
spectrum when $n \geq 2$.

%
%
\section{Algebras and modules}
\label{alg-mod}

A $\Gamma$-ring is a monoid in the category of $\Gamma$-spaces, i.e., a $\Gamma$-space $R$ equipped with a unit map 
$\s \to R$ and a multiplication map $R \wedge R \to R$ making the usual diagrams commute (see e.g. Mac Lane \cite[VII.3]{maclane}.) 
Given a $\Gamma$-ring $R$, we can consider the category of modules over $R$.  A left $R$-module is a $\Gamma$-space $M$
with an action $R \wedge M \to M$, again making certain obvious diagrams commute, and maps of $R$-modules are maps of
$\Gamma$-spaces that respect the action.  We let $\gmod_R$ denote the category of left $R$-modules.
Given a commutative $\Gamma$-ring $R$, we have the category of algebras over $R$.  An $R$-algebra is a monoid in the category of $R$-modules, and maps of $R$-algebras are maps of $R$-modules respecting the monoid structure.
Let $\galg_R$ denote the category of $R$-algebras.

Since $\gspc$ satisfies the monoid axiom, we can apply \cite[4.1]{schwede-shipley} and immediately get model structures
on the categories of modules and algebras over a monoid.
Here we are assuming the stable model structure on $\gspc$.  Of course, the result is also true for the 
pointwise model structure.

\begin{thm}
\label{alg-mod-model}
Suppose the model structure on $\spc$ is monoidal, and let $R$ be a $\Gamma$-ring.  Then the category  $\gmod_R$ inherits a cofibrantly generated model structure from $\gspc$.  If $R$ is commutative the same result 
holds for the category $\galg_R$, and every cofibrant $R$-algebra is also cofibrant as an $R$-module.
\end{thm}

The model structures in Theorem \ref{alg-mod-model} are created by forgetful functors: a map $f$ of $R$-modules is a 
weak equivalence (fibration) if and only if
its image $Uf$ under the forgetful functor $U \colon \gmod_R \to \gspc$ is a weak equivalence (fibration).  Similarly 
for $R$-algebras.

As an application we now establish some results about the Eilenberg-Mac Lane $\Gamma$-spaces, and the correspondence
with presheaves of simplicial abelian groups and rings.  The following are the presheaf versions of results in Schwede 
\cite{schwede}.
Let $\sabpre$ be the category of presheaves simplicial abelian groups.  For a monoid $A$ in $\sabpre$ let $\smodpre$ 
be the category of $A$-modules, and for a commutative monoid $B$ let $\salgpre$ be the category of $B$-algebras.
A map in $\sabpre$ is a weak equivalence (fibration) if the underlying map of spaces is a local weak equivalence 
(fibration.)  In the same way, weak equivalences and fibrations in $\smodpre$ and $\salgpre$ are defined on the 
underlying spaces.

\begin{thm}
With the above definitions of weak equivalences and fibrations, the category $\sabpre$ is a
cofibrantly generated model category, with generating cofibrations $\Z(I)$ and generating trivial cofibrations $\Z(J)$.
If $\spc$ is monoidal, the categories $\smodpre$ and $\salgpre$ are cofibrantly generated model categories as well.
\end{thm}
\begin{proof}
The category $\sabpre$ is bicomplete, and the rectract and 2-out-of-3 axiom follow immediately; we have to prove the second half
of the lifting axiom and the factorization axiom.  These follow by a standard argument involving (a transfinite version 
of) Quillen's small object argument.  

Let
$$
\Z : \spc \rightleftarrows \sabpre : U
$$
be the adjoint pair consisting of the free simplicial abelian presheaf functor $\Z$ and the forgetful functor $U$.
First note that maps in $\Z(I)$ are cofibrations in $\sabpre$, since by adjointness lifts in diagrams of the form
$$
\xymatrix{
\Z(A) \ar[r] \ar[d] & X \ar[d] \\
\Z(B) \ar[r] & Y
}
$$
are in a one-to-one correspondance with lifts in diagrams of the form
$$
\xymatrix{
A \ar[r] \ar[d] & U(X) \ar[d] \\
B \ar[r] & U(Y) \,.
}
$$
Also, maps in $\Z(J)$ are trivial cofibrations since by \cite[2.1]{jardine-chain} the functor $\Z$ preserves weak equivalences.

For the factorization axiom, let $f : X \to Y$ be a map in $\sabpre$.  By \cite[2.1.14]{hovey-model} the map $f$ 
can be factored as $f=p \circ i$, where $i$ is in $\Z(I)$-cell and
$p$ has the right lifting property with respect to maps in $\Z(I)$.  Since cofibrations in $\sabpre$ are defined by a left lifting
property, maps in $\Z(I)$-cell are cofibrations and in particular $i$ is a cofibration.
By adjointness $U(p)$ has the right lifting property with respect to maps in $I$, so $U(p)$ is a trivial fibration in $\spc$ and hence $p$
is a trivial fibration in $\sabpre$.  The other half of the factorization axiom is proved in a similar way, once
we know that maps in $\Z(J)$-cell are trivial cofibrations in $\sabpre$.  But all maps in $\Z(J)$ are 
monomorphisms and local weak equivalences, i.e., trivial cofibrations in $\spc$ with the local injective model
structure.  Trivial cofibrations are closed under forming cell objects, so in particular maps in $\Z(J)$-cell are 
local weak equivalences.

The last lifting axiom follows, since now we can factor each trivial cofibration $i$ as a map $j$ in $\Z(J)$-cell followed by a trivial fibration $p$.  There is a lift in the diagram
$$
\xymatrix{
A \ar[d]_i \ar[r]^j & C \ar[d]_p \\
B \ar[r] & B
}
$$
which shows that $i$ is a retract of $j$.  Maps in $\Z(J)$ have the left lifting property with respect to
fibrations, and this is also true for maps in $\Z(J)$-cell.  Since $i$ is a retract of $j$ we conclude that it has
the required lifting property.

The model structures for $\smodpre$ and $\salgpre$ follow from \cite[4.1]{schwede-shipley}.
\end{proof}

Let $A$ be a presheaf of simplicial abelian groups.  The Eilenberg-Mac Lane $\Gamma$-space $HA$ associated to $A$ is
defined as follows.  For each $n_+$ in $\Gamma$ let $HA(n_+) = A^{\times n}$, and for each map $f\colon n_+ \to m_+$
let the induced map $HA(n_+) \to HA(m_+)$ be defined by
$$
(a_1, \dots, a_n) \mapsto (\sum_{f(i)=1}a_i , \dots, \sum_{f(i)=m}a_i)
$$
in each section.  A map of simplicial abelian presheaves $A \to B$ induces a map of $\Gamma$-spaces $HA \to HB$.
Note that $HA$ is very special, and its associated spectrum is a generalized Eilenberg-Mac Lane spectrum for $A$ 
since $\pi_n(HA) = \pi_n(HA(1_+)) = \pi_n(A)$.

A functor $L$ in the opposite direction is described as follows.  Let $F$ be a $\Gamma$-space, and consider the
map
\begin{equation}
\label{L-cok}
p_{1*} + p_{2*} - \nabla_* \colon \widetilde\Z F(2_+) \to \widetilde\Z F(1_+) \,,
\end{equation}
where $p_1$ and $p_2$ are the two projections $2_+ \to 1_+$ in $\Gamma$, $\nabla$ is the fold map,
and $\widetilde\Z$ denotes the reduced free
simplicial abelian presheaf associated to a space.  The value of $L$ on $F$ is now defined to be the
cokernel of \ref{L-cok}.

The following result is just a sectionwise application of \cite[1.2]{schwede}.

\begin{lem}
\label{hl-adj}
The functor $L$ is strong symmetric monoidal, while $H$ is lax symmetric monoidal.
There is an adjunction
$$
L : \gspc \rightleftarrows \sabpre : H \,.
$$
Both $L$ and $H$ preserve modules, rings, and commutative rings.
Let $A$ be a presheaf of simplicial rings and $B$ be a presheaf of commutative simplicial rings.  The functors $L$ and
$H$ induces adjunctions
\begin{eqnarray*}
L : \gmod_{HA} \rightleftarrows \smodpre : H \\
L : \galg_{HB} \rightleftarrows \salgpre : H \,.
\end{eqnarray*}
\end{lem}

\begin{lem}
All three adjunctions in Lemma \ref{hl-adj} are Quillen adjunctions.
\end{lem}
\begin{proof}
Let us consider the first adjunction, the result for the other two follows by the same argument.
Since trivial fibrations of spaces are closed under finite products, $H$ takes trivial fibrations of simplicial
abelian presheaves to pointwise trivial fibrations of $\Gamma$-spaces, which coincide with the stably trivial 
fibrations of $\Gamma$-spaces.

The functor $H$ also takes fibrations of simplicial abelian presheaves to pointwise fibrations of $\Gamma$-spaces
between stably fibrant $\Gamma$-spaces, which coincides with stable fibrations between stably fibrant $\Gamma$-spaces.
\end{proof}

\begin{thm}
\label{simp-mod}
Let $A$ be a presheaf of simplical rings.  Then the adjoint functors $H$ and $L$ constitute a Quillen equivalence 
between the categories of presheaves of simplicial $A$-modules and $HA$-modules.
\end{thm}
\begin{proof}
The following proof is an adaption of Schwede's argument given in \cite[4.2]{schwede}.
The functor $H$ preserves weak equivalences, and detects weak equivalences since a stable equivalence $HM \to HN$ is
a pointwise equivalence, and in particular $M=HM(1_+) \to HN(1_+)=N$ is a local weak equivalence.  It remains to show
that for every cofibrant $HA$-module $M$ the unit map $M \to HL(M)$ is a stable equivalence.

We first consider $\Gamma$-spaces of the form $HA \wedge X$, where $X$ is a space, and we claim that the presheaf map
$\pi^p_*(HA \wedge X) \to \pi^p_*(HL(HA \wedge X))$ is a sectionwise isomorphism.  After evaluating in sections
we are led to consider the map $\pi_*(HA(U) \wedge K) \to \pi_*(HL(HA \wedge K)(U))$ as a natural transformation of
functors of the simplicial set $K$.  But this is easily seen to be an isomorphism for the case $K=S^0$, and both
functors are homology theories with coefficients in $A$, since $L(HA \wedge K)(U)$ is just the free $A(U)$-module
generated by $K$.  Thus the map is an isomorphism for all $K$ and in particular for $X(U)$.

The map $\Gamma^1 \wedge n_+ \to \Gamma^n$ induced by the $n$ projections $n_+ \to 1_+$ is a stable equivalence,
since the induced map of spectra is just the canonical inclusion $\s^{\vee n} \to \s^{\times n}$.  This implies that
$F \wedge n_+ \cong F \wedge \Gamma^1 \wedge n_+$ is stably equivalent to $F \wedge \Gamma^n$ for 
all $\Gamma$-spaces $F$.  The composite functor $HL$ preserves weak equivalences between cofibrant objects, so the
unit map of $HA \wedge X \wedge \Gamma^n$ is a stable equivalence by the case already proved.

Let $M$ be a cofibrant $HA$-module, i.e., a retract of a colimit $\colim_{\alpha<\gamma} M_\alpha$, where $\gamma$ 
is an ordinal and the maps $M_\alpha \to M_{\alpha+1}$ are pushouts of generating cofibrations in $\gmod_{HA}$.  
The generating cofibrations in $\gmod_{HA}$ are of the form 
$$
HA \wedge X \wedge \Gamma^n \to HA \wedge Y \wedge \Gamma^n \,,
$$
where $X \to Y$ is a (generating) cofibration of spaces.  If we have a pushout diagram of the form
$$
\xymatrix{
HA \wedge X \wedge \Gamma^n \ar[r] \ar[d] & M_\alpha \ar[d] \\
HA \wedge Y \wedge \Gamma^n \ar[r] & M_{\alpha+1}
}
$$
and assume that the map $M_\alpha \to HL(M_\alpha)$ is a stable equivalence, we can use the first part and the gluing
lemma (see e.g. \cite[II.8.12]{gj}) to show that the map $M_{\alpha+1} \to HL(M_{\alpha+1})$ is a stable equivalence.
Now the induced map 
$$
\colim_{\alpha<\gamma} M_\alpha \to \colim_{\alpha<\gamma} HL(M_\alpha)
$$ 
is a stable equivalence, and $\colim HL(M_\alpha)$ is stably equivalent to $HL(\colim M_\alpha)$ since $L$ preserves
colimits and 
$$
\pi_*(\colim HA_\alpha) \cong \colim\pi_*(HA_\alpha) \cong \colim \pi_*(A_\alpha) \cong \pi_*(\colim A_\alpha) 
\cong \pi_*H(\colim A_\alpha) \,.
$$
Finally, since $M$ is a retract of $\colim M_\alpha$, the unit map $M \to HL(M)$ is also a stable equivalence.
\end{proof}

\begin{thm}
Let $B$ be a presheaf of commutative simplical rings.  Then the adjoint functors $H$ and $L$ are a Quillen equivalence 
between the categories of presheaves of simplicial $B$-algebras and $HB$-algebras.
\end{thm}
\begin{proof}
Since every cofibrant $HB$-algebra is cofibrant as an $HB$-module, the proof of Theorem \ref{simp-mod} applies.
\end{proof}

\numberwithin{equation}{subsection}

%
%
\section{Appendix}

\subsection{Simplicial spaces}

Given a simplicial space $X$, i.e., a bisimplicial presheaf, we obtain a space $X_{m,*}$ by fixing the first 
simplicial degree $m$.  We say that a map $X \to Y$ is a pointwise equivalence if $X_{m,*} \to Y_{m,*}$ is a 
local weak equivalence for all $m$.

\begin{prop}
\label{diagonal-eq}
Let $X \to Y$ be a pointwise equivalence of simplicial spaces.  Then the induced diagonal map $dX \to dY$ is a local 
weak equivalence of spaces.
\end{prop}
\begin{proof}
The result only depends
on the weak equivalences on simplicial presheaves, so we are free to choose the local injective model structure
where every object
is cofibrant.  Now the proof in \cite[IV.1.7]{gj} for bisimplicial sets carries over, mutatis mutandis.  
\end{proof}

\subsection{Enriched left Bousfield localization}

Here we summarize the theory of enriched left Bousfield localization as developed in Barwick \cite{barwick}.
We will ignore the set-theoretic details that appear in these statements; they are treated carefully in Barwick's paper.

\begin{defn}
\label{def-bous}
Let $\vcat$ be a monoidal model category and $\cat$ a $\vcat$-model category.  Suppose $\Sigma$ is a set of
morphisms in $\cat$.  A left Bousfield localization of $\cat$ with respect to $\Sigma$ enriched over $\vcat$ is
a $\vcat$-model category $L_{\Sigma/\vcat} \cat$, equipped with a left Quillen $\vcat$-functor $\cat \to L_{\Sigma/\vcat} \cat$ 
that is initial among left Quillen $\vcat$-functors $L \colon \cat \to \mathcal N$ to $\vcat$-model categories $\mathcal N$ such
that $Lf$ is a weak equivalence in $\mathcal N$ for all $f$ in $\Sigma$.
\end{defn}

\begin{defn}
Let $\vcat$, $\cat$ and $\Sigma$ be as in Definition \ref{def-bous}.
\begin{itemize}
\item An object $Z$ in $\cat$ is \emph{$\Sigma/\vcat$-local} if it is fibrant, and for any morphism $A \to B$ in $\Sigma$ the morphism
$$
\vhom(B_c,Z) \to \vhom(A_c,Z)
$$
is a weak equivalence in $\vcat$.
\item A morphism $A \to B$ in $\cat$ is a \emph{$\Sigma/\vcat$-local equivalence} if for any  
$\Sigma/\vcat$-local object $Z$ in $\cat$, the morphism
$$
\vhom(B_c,Z) \to \vhom(A_c,Z)
$$
is a weak equivalence in $\vcat$.
\end{itemize}
\end{defn}

The following result is proved in \cite[3.18]{barwick}.
\begin{thm}
\label{enriched-bousfield}
Suppose that $\vcat$ is a combinatorial monoidal model category and $\cat$ is a left proper and combinatorial 
$\vcat$-model category.  Suppose further that the generating cofibrations and generating trivial cofibrations in $\vcat$ 
and $\cat$ all have cofibrant domains.  Let $\Sigma$ be a set of morphisms in $\cat$.
Then the left Bousfield localization of $\cat$ with respect to $\Sigma$ enriched over $\vcat$ exists, and it has the
following properties.
\begin{itemize}
\item As a category, $L_{\Sigma/\vcat} \cat$ is just $\cat$.
\item The model category $L_{\Sigma/\vcat} \cat$ is combinatorial and left proper.
\item The cofibrations in $L_{\Sigma/\vcat} \cat$ are the same as those of $\cat$.
\item The fibrant objects in $L_{\Sigma/\vcat} \cat$ are the fibrant $\Sigma/\vcat$-local objects in $\cat$.
\item The weak equivalences in $L_{\Sigma/\vcat} \cat$ are the $\Sigma/\vcat$-local equivalences.
\end{itemize}
\end{thm}


\bibliographystyle{gtart}
\bibliography{gamma-gt}

\end{document}